\documentclass[12pt]{amsart}
\renewcommand{\oe}{\operatorname{\mbox{\rm{\OE}}}}
\newtheorem{thm}{Theorem}
\newtheorem{cor}{Corollary}
\newtheorem{prop}{Proposition}
\begin{document}
\title[Smoothly Parameterised \v{C}ech Cohomology]
{Smoothly Parameterised \v{C}ech Cohomology\\ of Complex Manifolds}
\author{Toby Bailey}
\address{\hskip-\parindent
Department of Mathematics\\
University of Edinburgh\\
James Clerk Maxwell Building\\
The King's Buildings\\
Mayfield Road\\
Edinburgh EH9 3JZ\\
Scotland}
\email{tnb@mathematics.edinburgh.ac.uk}
\author{Michael Eastwood}
\address{\hskip-\parindent
Pure Mathematics Department\\
Adelaide University\\
South Australia 5005}
\email{meastwoo@maths.adelaide.edu.au}
\author{Simon Gindikin}
\address{\hskip-\parindent
Department of Mathematics\\
Rutgers University\\
New Brunswick\\
NJ 08903\\
USA}
\email{gindikin@math.rutgers.edu}
\renewcommand{\subjclassname}{\textup{2000} Mathematics Subject Classification}
\renewcommand{\subjclassname}{\textup{2000} Mathematics Subject Classification}
\thanks{Support from the Australian Research Council and the Mathematical
Sciences Research Institute is gratefully acknowledged.
Research at MSRI is supported in part by NSF grant DMS-9810361.
Simon Gindikin is partially supported by NSF grant DMS-0070816.}
\begin{abstract}
A Stein covering of a complex manifold may be used to realise its analytic
cohomology in accordance with the \v{C}ech theory. If, however, the Stein
covering is parameterised by a smooth manifold rather than just a discrete set,
then we construct a cohomology theory in which an exterior derivative replaces
the usual combinatorial \v{C}ech differential. Our construction is motivated by
integral geometry and the representation theory of Lie groups.
\end{abstract}
\maketitle
\section{Introduction}
The usual language of \v{C}ech cohomology is adapted for discrete coverings
$\{U_i\}_{i\in I}$: the parameter space $I$ has no particular structure. In
complex analysis, however, it is typical to consider infinite coverings
$\{U_\xi\}_{\xi\in\Xi}$ of a complex manifold $Z$ by open Stein subsets which
are themselves parameterised by points of an auxiliary smooth manifold~$\Xi$.
In such a setting, it is unnatural to forget the nature of the parameter space
and treat the parameters as discrete. Instead, it is natural to consider
de~Rham cohomology on $\Xi$ depending holomorphically on points in elements of
the covering. The exact details are in the following section. It is the
smoothly parameterised \v{C}ech cohomology introduced in~\cite{g2} and further
discussed in~\cite{g3}. In suitable very general circumstances, we obtain the
analytic cohomology $H^p(Z,{\mathcal O})$ of~$Z$. This transfer from the usual
combinatorial \v{C}ech language to de~Rham cohomology on the parameter space
$\Xi$ with holomorphic coefficients in $U_\xi$ is akin to the transition from
integral sums to integrals.

Smoothly parameterised \v{C}ech cohomology provides an alternative to the
classical \v{C}ech and Dolbeault realisations of $H^p(Z,{\mathcal O})$.
In fact, not only is it one more natural way to present analytic cohomology
but, as we will see in examples, it is a very effective way to connect analytic
cohomology with the geometry of complex manifolds.

Let us emphasise that using continuous \v{C}ech cohomology instead of discrete,
we have a very explicit way to construct an operator from \v{C}ech cohomology to
Dolbeault cohomology using a smooth `section', viz.\ $\gamma:Z\to\Xi$ such that
$z\in U_{\gamma(z)},\,\forall z\in Z$. This direct link with Dolbeault
cohomology is unavailable in usual \v{C}ech theory.

The first aim of this article is to make this informal discussion completely
precise and to state the technical conditions under which the $p^{\mathrm{th}}$
cohomology of the resulting `de~Rham complex with holomorphic coefficients' is
$H^p(Z,{\mathcal O})$. In fact, the proof extends to a considerably more
general setting. In effect, the covering need not be by open subsets but,
rather, by Stein submanifolds. At the same time, the parameter space $\Xi$ is
not strictly necessary: the appropriate submersion is replaced by a foliation.
Our general formulation also deals with the case when the parameter space $\Xi$
is itself a complex manifold: we obtain the theory established
in~\cite{egw,egw1}. The main analytic input for our general proof is due to
Jurchescu~\cite{j}. Other ingredients for our construction are taken
from~\cite{bes}. The results of this article were roughly sketched
in~\cite{beg}---here we provide complete details.

We state our basic result and prove it in a more general setting than for
smoothly parameterised \v{C}ech cohomology. There are two reasons. Firstly, the
general proof is more natural. Secondly, we expect the general machinery to
apply to examples arising naturally in representation theory where there is no
smoothly parameterised Stein covering. In the interests of readers only
concerned with the smoothly parameterised \v{C}ech theory, we carefully define
this basic theory in~\S\ref{simple} and outline its proof in~\S\ref{summary}.
The full proof is a consequence of the general formulation and reasoning
in~\S\ref{general}.

Our primary motivation is the study of cohomology of tubes over non-convex
cones and the explicit representations this cohomology provides. We provide
some examples of our construction in this setting arising from pseudo-Hermitian
symmetric spaces.

\section{Formulation of the smooth \v{C}ech theory}\label{simple}
This section contains a formulation of smoothly parameterised \v{C}ech
cohomology and a statement of our main result, Theorem~\ref{SmoothCechTheorem}.
It provides sufficient conditions that the smoothly parameterised \v{C}ech
theory compute the usual analytic cohomology.

Let $Z$ be a complex manifold whose cohomology $H^p(Z,{\mathcal O})$ we wish to
compute. Suppose we are given an open covering $\{U_\xi\}_{\xi\in\Xi}$ of $Z$
where $\Xi$ is a smooth manifold. We need a sense in which $U_\xi$ depends
smoothly on~$\xi$. For this, let us introduce
$$M=\{(\xi,z)\in\Xi\times Z\mbox{ s.t. }z\in U_\xi\}\subseteq\Xi\times Z$$
and simply insist that $M$ be an open subset of $\Xi\times Z$.

\subsection{First example}\label{ichi} A useful example to bear in mind is
$Z={\mathbb C}^n\setminus{\mathbb R}^n$. In this case, we may take
\begin{equation}\label{one}
U_\xi=\{z=x+iy\in{\mathbb C}^n\mbox{ s.t. }\langle\xi,y\rangle>0\},
\quad\mbox{for }
\xi\in\Xi=\{\xi\in{\mathbb R}^n\mbox{ s.t. }|\xi|=1\}
\end{equation}
as an especially natural covering. Notice that each $U_\xi$ is convex and hence
Stein. Also, for each $z\in Z$ the set $\{\xi\in\Xi\mbox{ s.t. }z\in U_\xi\}$
is a hemisphere. Thus, $M$ is an open subset of $S^{n-1}\times Z$ and
the natural projection $M\to Z$ has contractible fibres. A link between the
smoothly parameterised \v{C}ech theory in this case and Sato's cohomological
description of hyperfunctions on ${\mathbb{R}}^n$ is provided in~\cite{beg,g3}.

\subsection{Second example}\label{ni} Let $Z={\mathbb{CP}}_n\setminus B$, the
complement of a closed ball in complex projective space:--
$$Z=\{[z_0,z_1,\ldots,z_n]\in{\mathbb{CP}}_n\mbox{ s.t. }
|z_0|^2<|z_1|^2+\cdots+|z_n|^2\}.$$
Consider $\xi\in{\mathrm{Gr}}_{n-1}({\mathbb{C}}^{n+1})$ as a linearly embedded
${\mathbb{CP}}_{n-2}\hookrightarrow{\mathbb{CP}}_n$ and let $\Xi$ denote those
that avoid~$B$. Then, for each $\xi\in\Xi$ we may take
$$U_\xi=\{z\in Z
\mbox{ s.t. the join of $z$ and $\xi$ continues to avoid $B$}\}.$$
This is a smoothly parameterised Stein covering of~$Z$. Even better, $\Xi$ is
itself a complex manifold and the holomorphic language of \cite{egw,egw1,g1}
may be employed. This language may be employed easily to invert the
Penrose-Radon transform (the case $n=2$ is detailed in~\cite{beg}).

\subsection{Third example}\label{san} Let $V\subset{\mathbb{R}}^3$ be the
following non-convex cone:--
$$V=\{(x_1,x_2,x_3)\in{\mathbb{R}}^3\mbox{ s.t. }x_1{}^2+x_2{}^2>x_3{}^2\}.$$
It is one of the open orbits of the standard action of ${\mathrm{SO}}(2,1)$
on ${\mathbb{R}}^3$ (the other two being convex cones). Let
$Z=V+i{\mathbb{R}}^3\subset{\mathbb{C}}^3$ be the corresponding tube domain.
Let $\Xi$ be the unit circle and for $\theta\in\Xi$ take
$$U_\theta=\{z=x+iy\in{\mathbb{C}}^3\mbox{ s.t. }
x_1\cos\theta+x_2\sin\theta>|x_3|\}.$$
It is a tube over a convex cone. Thus, we have a Stein covering smoothly
parameterised by~$\Xi$. Notice, however, that this covering is not preserved by
the action of ${\mathrm{SO}}(2,1)$. The Stein covering
$$U_{\theta,\phi}=\{x+iy\mbox{ s.t. }x_1\cos\theta+x_2\sin\theta>x_3
\mbox{ and }x_1\cos\phi+x_2\sin\phi>-x_3\}$$
for $\theta\not=\phi$ remedies this.

\medbreak

These examples illustrate various points. In \S\ref{ichi} a finite Stein
covering by half spaces is available but cannot respect the evident rotational
symmetry of~$Z$. The smoothly parameterised covering is certainly more natural
from this point of view. The holomorphically parameterised covering in
\S\ref{ni} is perhaps the best choice. But, for the complex manifolds in
\S\ref{ichi} and~\S\ref{san}, the smoothly parameterised coverings
seem to be optimal. Thus, we are obliged to extend the holomorphic machinery of
\cite{egw,egw1} to the smooth case.

Returning to the general discussion of the smoothly parameterised case, it is
convenient to introduce submersions $\eta$ and~$\tau$
\begin{equation}\label{doublefibration}\raisebox{-25pt}{\begin{picture}(80,50)
\put(0,5){\makebox(0,0){$Z$}} \put(40,45){\makebox(0,0){$M$}}
\put(80,5){\makebox(0,0){$\Xi$}} \put(30,35){\vector(-1,-1){20}}
\put(50,35){\vector(1,-1){20}} \put(15,30){\makebox(0,0){$\eta$}}
\put(65,30){\makebox(0,0){$\tau$}}
\end{picture}}\end{equation}
as restrictions of the natural projections to~$M$. Locally, on $\Xi\times Z$ we
may consider complex-valued smooth functions $f=f(\xi,z)$ that are holomorphic
in $z$ for fixed~$\xi$. We shall say that such functions are {\em partially
holomorphic\/} and write $\oe$ for the sheaf of germs thereof. On $\Xi\times Z$
and hence also on~$M$, we may consider $p$-forms in the $\xi$-variables alone:
precisely, if we let $\Lambda_\Xi^p$ denote the bundles of $p$-forms on $\Xi$
and set $B^p=\tau^*\Lambda_\Xi^p$, then such $p$-forms are smooth sections
of~$B^p$. Geometrically, these are $p$-forms along the fibres of~$\eta$. We
shall also refer to them as $p$-forms {\em relative\/} to~$\eta$, or simply
relative $p$-forms if $\eta$ is understood. In smooth local
co\"ordinates $(\xi_1,\ldots,\xi_m)$ on $\Xi$ and holomorphic local
co\"ordinates $(z_1,\ldots,z_n)$ on $Z$, relative $p$-forms may be written
$$\sum_{i_1,\ldots,i_p}
\omega_{i_1\cdots i_p}(\xi_1,\ldots,\xi_m,z_1,\ldots,z_n)\,
d\xi_{i_1}\wedge\cdots\wedge d\xi_{i_p}$$
but it makes co\"ordinate-free sense to require that all the coefficients be
holomorphic in $z_1,\ldots,z_n$. We shall refer to such relative $p$-forms as
partially holomorphic and write $\oe(B^p)$ for the sheaf of germs thereof. It
is easily verified that the exterior derivative in the $\xi$-variables alone,
given by the usual formula with $z\in Z$ as a passenger, is well-defined and
takes partially holomorphic relative $p$-forms to partially holomorphic
relative $(p+1)$-forms. This is called the relative exterior derivative and is
written~$d_\eta$. Of course, $d_\eta{}^2=0$. In summary, we have a complex of
sheaves on~$M$
$$0\to\oe\stackrel{d_\eta}{\longrightarrow}\oe(B^1)
\stackrel{d_\eta}{\longrightarrow}\oe(B^2)\to\cdots\to
\oe(B^p)\stackrel{d_\eta}{\longrightarrow}\oe(B^{p+1})\to\cdots$$
and we may define the $p^{\mathrm{th}}$ smooth \v{C}ech cohomology to be
$$\frac{\ker:\Gamma(M,\oe(B^p))
             \stackrel{d_\eta}{\longrightarrow}\Gamma(M,\oe(B^{p+1}))}
       {\mbox{im : }\Gamma(M,\oe(B^{p-1}))
             \stackrel{d_\eta}{\longrightarrow}\Gamma(M,\oe(B^{p}))}.$$
In order that this coincide with the analytic cohomology $H^p(Z,{\mathcal O})$,
we need a smooth analogue of the Leray condition from usual (discrete) \v{C}ech
cohomology. Roughly speaking, we would like each $U_\xi$ to be Stein but in a
way that depends smoothly on $\xi\in\Xi$. To make this precise, we define
the {\em partially holomorphic hull\/} $\widehat K_M$ of a compact subset
$K\subset M$ to be
\begin{equation}\label{hull}\widehat K_M=\{x\in M\mbox{ s.t. }
|f(x)|\leq\sup_K|f|\;\;\forall f\in\Gamma(M,\oe)\}.\end{equation}
Then $M$ is a {\em Cartan manifold\/} in the sense of Jurchescu~\cite{j} if and
only if the following three conditions hold.
\begin{itemize}
    \item if $K\subset M$ is compact, then so is $\widehat K_M$;
\item the partially holomorphic functions on $M$ separate points;
\item the partially holomorphic functions on $M$ provide local co\"ordinates.
\end{itemize}
Specifically, local co\"ordinates $(\xi_1,\ldots,\xi_m,z_z\ldots,z_n)$ are each
partially holomorphic and the third requirement of a Cartan manifold is that
charts may be chosen in which the co\"ordinate functions extend to global
partially holomorphic functions. If $M$ is an open subset of
$\Xi\times{\mathbb C}^n$, then the final two conditions are evidently
satisfied. In this case, a Cartan manifold may reasonably be called a `domain
of partial holomorphy'.

For a smoothly parameterised Leray theorem, some topological restriction is
also necessary. For each $z\in Z$, the relative de~Rham sequence is trying to
compute the de~Rham cohomology of the fibre~$\eta^{-1}(z)$. We need to
eliminate this effect:--
\begin{thm}\label{SmoothCechTheorem}
Let $Z$ be a complex manifold endowed with a smoothly parameterised open cover
viewed as~\mbox{\rm (\ref{doublefibration})}. Suppose $M$ is a Cartan manifold
and $\eta$ has contractible fibres. Then
$$H^p(Z,{\mathcal O})\cong H^p(\Gamma(M,\oe(B^\bullet))),\ \forall p.$$
\end{thm}
\noindent This theorem is a special case of a much more general result,
Theorem~\ref{mainthm}, whose formulation and proof will occupy the next two
sections.

In order for Theorem~\ref{SmoothCechTheorem} to be useful, of course, we need a
good supply of Cartan manifolds:--
\begin{prop}\label{supply}
Let $\Xi$ be a smooth manifold. If $M\subset\Xi\times{\mathbb C}^n$
is open and the natural projection $\tau:M\to\Xi$ has convex fibres, then $M$
is Cartan.
\end{prop}
\begin{proof}It suffices to check partial holomorphic convexity: the other two
requirements are inherited from $\Xi\times{\mathbb C}^n$. For each linear
functional $\ell:{\mathbb C}^n\to{\mathbb C}$ and $\phi(\xi)$ a smooth bump
function on~$\Xi$, the function
$$(\xi,z)\mapsto\phi(\xi)\exp\ell(z)$$
is partially holomorphic. Hence, the partially holomorphic hull of any subset
$K\subset M$ is contained in its fibrewise convex hull.
\end{proof}
Notice that if $M\subset\Xi\times{\mathbb C}^n$ is a Cartan manifold, then the
fibres of $\tau:M\to\Xi$ are Stein. The converse, however, is false:--
$$M={\mathbb R}\times{\mathbb C}\setminus\{(0,0)\}
\stackrel{\tau}{\longrightarrow}{\mathbb R}$$
has Stein fibres but is not Cartan. The Cauchy integral formula shows that
partially holomorphic functions extend across the origin. Consequently,
$$K=\{(\xi,z)\in M\mbox{ s.t. }|\xi|=1,\;|z|\leq 1\}\implies
\widehat K_M=\{(\xi,z)\in M\mbox{ s.t. }|\xi|\leq 1,\;|z|\leq 1\},$$
which is not compact.

Suppose the hypotheses of Theorem~\ref{SmoothCechTheorem} are satisfied and we
are given $\gamma:Z\to M$ a smooth section of the submersion $\eta:M\to Z$. We
may restrict forms on $Z$ to the image of this section, consider the result as
a form on~$Z$, and take the $(0,p)$-part. This procedure gives a more explicit
interpretation of Theorem~\ref{SmoothCechTheorem} in terms of Dolbeault
cohomology. Formally, we may proceed as follows. A smooth section $\omega$ of
$B^p$ is, in particular, a $p$-form on~$M$ whence its pullback $\gamma^*\omega$
is a $p$-form on~$Z$. Let ${\mathcal E}_Z^{0,p}$ denote the sheaf of smooth
forms on $Z$ of type~$(0,p)$ and define
\begin{equation}\label{pullback}
\Gamma(M,\oe(B^p))\longrightarrow\Gamma(Z,{\mathcal E}_Z^{0,p})
\end{equation}
by $\omega\mapsto(\gamma^*\omega)^{0,p}$, the $(0,p)$-component
of~$\gamma^*\omega$. As $M$ is locally a product, we can split the exterior
derivative into vertical and horizontal parts:--
$$d\omega=d_\eta\omega+d_\tau\omega=
          d_\eta\omega+\partial_\tau\omega+\bar\partial_\tau\omega,$$
where $d_\tau$ is further split according to type. That $\omega$ is partially
holomorphic is to say $\bar\partial_\tau\omega=0$. Evidently,
$\gamma^*\partial_\tau\omega$ is of type $(1,p)$ on~$Z$. Therefore,
$$(\gamma^*d_\eta\omega)^{(0,p+1)}=(\gamma^*d\omega)^{(0,p+1)}=
(d(\gamma^*\omega))^{(0,p+1)}=\bar\partial((\gamma^*\omega)^{(0,p)}).$$
In other words, (\ref{pullback}) is a homomorphism of complexes. The right hand
side is the Dolbeault complex on $Z$. The following result is a special case of
Theorem~\ref{sectionthm} in~\S\ref{general}.
\begin{thm} If the hypotheses of Theorem~\ref{SmoothCechTheorem} are satisfied
and $\gamma:Z\to M$ is an arbitrary smooth section of the submersion
$\eta:M\to Z$, then the homomorphism {\rm{(\ref{pullback})}} induces an
isomorphism on cohomology.
\end{thm}

To close this discussion, let us revisit our three examples. In \S\ref{ichi}
the covering (\ref{one}) of ${\mathbb C}^n\setminus{\mathbb R}^n$ satisfies the
conditions of Theorem~\ref{SmoothCechTheorem} and admits a natural smooth
section~$\gamma$. We have observed already that the fibres of $\eta$ are, as
hemispheres, contractible. That
$$M=\{(\xi,z=x+iy)\in{\mathbb R}^n\times{\mathbb C}^n\mbox{ s.t. }
|\xi|=1\mbox{ and }\langle\xi,y\rangle>0\}$$
is a Cartan manifold follows immediately from Proposition~\ref{supply}. For
$\gamma:{\mathbb C}^n\setminus{\mathbb R}^n\to M$ we may take
$z=x+iy\mapsto\xi=y/|y|$.

In \S\ref{ni} the corresponding manifold $M$ is Stein, as it should be in the
case of a holomorphically parameterised covering. For $\gamma:Z\to M$ we may
take
$$z=[z_0,z_1,\ldots,z_n]\stackrel{\gamma}{\longmapsto}
(\{[0,\xi_1,\ldots,\xi_n]\in{\mathbb{CP}}_n\mbox{ s.t. }
\bar{z_1}\xi_1+\cdots+\bar{z_n}\xi_n=0\},z)$$
and the fibres of $\eta$ are easily seen to contract onto this section. Note,
however, that $\gamma$ is only ${\mathrm{SU}}(n)$-invariant, whereas the full
symmetry group is ${\mathrm{SU}}(1,n)$.

In \S\ref{san} with either of the two coverings, the relevant manifold $M$ is
Cartan from Proposition~\ref{supply} and for $\gamma:Z\to M$ we may take
$$\gamma(x+iy)=(\arctan(x_2/x_1),x+iy)\quad\mbox{or}\quad
(\arctan(x_2/x_1),\arctan(x_2/x_1),x+iy),$$
easily checking that the fibres of $\eta$ contract onto its image. This section
$\gamma$ is only ${\mathrm{SO}}(2)$-invariant, whilst the full symmetry group
for the second covering is ${\mathrm{SO}}(2,1)$.

\section{Mixed manifolds}\label{mixed}
The reader only interested in smoothly parameterised \v{C}ech cohomology and
its applications will find a summary of the proof for this case
in~\S\ref{summary}.

A {\em mixed manifold} $M$ of type $(m,n)$ is a smooth manifold of dimension
$m+2n$ equipped with a Levi flat CR structure of codimension~$m$. More
precisely, we are given an integrable distribution $H\subset TM$ of rank $2n$
and an endomorphism $J:H\to H$ whose Nijenhuis tensor,
$$N(X,Y)\equiv[X,Y]+J[JX,Y]+J[X,JY]-[JX,JY]\quad
\mbox{for }X,Y\in\Gamma(M,H),$$
vanishes. Thus, $M$ is equipped with a foliation of dimension $2n$ with
smoothly varying complex structure on the leaves. If we set
$$T_M^{0,1}=\{X\in{\mathbb C}H\mbox{ s.t. }(J+i)X=0\},$$
then
$$T_M^{0,1}\oplus\overline{T_M^{0,1}}={\mathbb C}H\quad\mbox{and}\quad
[T_M^{0,1},T_M^{0,1}]= T_M^{0,1},$$
where this second statement means that $T_M^{0,1}$ is closed under Lie bracket.
Equivalently, we may set $\Lambda_M^{1,0}=(T_M^{0,1})^\perp$ and ask that this
generate a differentially closed ideal.

By the Newlander-Nirenberg theorem~\cite{nn}, mixed manifolds are locally
modelled on ${\mathbb R}^m\times{\mathbb C}^n$ with transition functions of the
form
$$(\xi,z)\mapsto(s(\xi),w(\xi,z)),$$
where $w(\xi,z)$ is holomorphic in~$z$. It is evident that mixed manifolds
provide a natural generalisation of the open subsets $M\subset\Xi\times Z$ that
we encountered in the previous section. Various notions extend immediately.
{\em Partially holomorphic\/} functions on a general mixed manifold $M$ are
smooth functions whose restriction to the leaves of the foliation are
holomorphic. In other words, these are the smooth complex-valued functions
annihilated by vector fields from~$T_M^{0,1}$. As before, we shall write $\oe$
for the sheaf of germs of partially holomorphic functions. If we define the
bundle of $(0,1)$-forms on $M$ by
$\Lambda_M^{0,1}:=\Lambda_M^1/\Lambda_M^{1,0}$ and write ${\mathcal E}_M^{0,1}$
for the corresponding sheaf of smooth sections, then we have a resolution
\begin{equation}\label{partialdolbeault}
0\to\oe\to{\mathcal E}_M^{0,0}\stackrel{\bar\partial}{\longrightarrow}
{\mathcal E}_M^{0,1}\stackrel{\bar\partial}{\longrightarrow}
{\mathcal E}_M^{0,2}\stackrel{\bar\partial}{\longrightarrow}\cdots
\stackrel{\bar\partial}{\longrightarrow}{\mathcal E}_M^{0,n}\to 0,
\end{equation}
which is the Dolbeault resolution on each leaf.

A mixed manifold of type $(0,n)$ is simply a complex manifold, in which case
(\ref{partialdolbeault}) is the Dolbeault resolution. The {\em partially
holomorphic hull\/} $\widehat K_M$ of a compact subset $K\subset M$ is defined
exactly as before by~(\ref{hull}) and the definition of a {\em Cartan
manifold\/} also reads the same. By analogy with the classical Levi problem on
a Stein manifold, Jurchescu~\cite{j} proved the following vanishing theorem:--
\begin{thm}\label{vanishingthm}If $M$ is a Cartan manifold, then
$$H^q(M,\oe)=0,\ \forall q\geq 1.$$
\end{thm}
\noindent An alternative proof may be found in~\cite{a}. We are grateful to
Gennadi Henkin for drawing our attention to these articles.

\section{General formulation and proof}\label{general}
Suppose $M$ is a mixed manifold, $Z$ is a complex manifold, and $\eta:M\to Z$
is a partially holomorphic submersion. More precisely, we shall suppose that
$\eta$ is a smooth surjection of maximal rank and that $\eta$ is holomorphic on
the leaves of the mixed manifold structure. We do not assume any particular
relationship between the dimension of $Z$ and the dimension of the complex
leaves in~$M$. Extreme special cases are when $M$ is itself complex or when $M$
is smooth (in which case it may be simply $Z$ but regarded as a smooth
manifold).

That $\eta$ is maximal rank is to say that $\eta^*\Lambda_Z^1\to\Lambda_M^1$ is
an injection.
That $\eta$ is partially holomorphic is to say that
$\eta^*\Lambda_Z^{1,0}\subseteq\Lambda_M^{1,0}$. Define a vector bundle $B^1$
on $M$ by the exact sequence
$$0\to\eta^*\Lambda_Z^{1,0}\to\Lambda_M^{1,0}\to B^1\to 0.$$
We shall see, in proving Theorem~\ref{mainthm} below, that $B^1$ is
naturally holomorphic on the leaves of the foliation. The same is true for
$B^p\equiv\Lambda^p(B^1)$. We shall say that such bundles are partially
holomorphic and write $\oe(B^p)$ for the sheaf of smooth sections of $B^p$ that
are holomorphic along the leaves. The following theorem generalises
Theorem~\ref{SmoothCechTheorem}.
\begin{thm}\label{mainthm} There is a complex
of sheaves on~$M$
\begin{equation}\label{Bcomplex}0\to\oe(B^0)
\stackrel{d_\eta}{\longrightarrow}\oe(B^1)
\stackrel{d_\eta}{\longrightarrow}\oe(B^2)\to\cdots\to
\oe(B^p)\stackrel{d_\eta}{\longrightarrow}\oe(B^{p+1})\to\cdots\end{equation}
so that if $M$ is a Cartan manifold and the fibres of $\eta$ are contractible,
then
\begin{equation}\label{mainresult}
H^p(Z,{\mathcal O})\cong H^p(\Gamma(M,\oe(B^\bullet))),\ \forall p.
\end{equation}
\end{thm}
\begin{proof}
In order to understand the bundle $B^1$ on~$M$, let us consider the following
commutative diagram with exact rows and columns.
\begin{equation}\label{diagram}\begin{array}{ccccccccc}
&&0&&0\\
&&\uparrow&&\uparrow\\
&&\eta^*\Lambda_Z^{0,1}&\to&\Lambda_M^{0,1}\\
&&\uparrow&&\uparrow\\
0&\to&\eta^*\Lambda_Z^1&\to&\Lambda_M^1&\to&\Lambda_\eta^1&\to&0\\
&&\uparrow&&\uparrow&&\uparrow\\
0&\to&\eta^*\Lambda_Z^{1,0}&\to&\Lambda_M^{1,0}&\to&B^1&\to&0\\
&&\uparrow&&\uparrow\\
&&0&&0
\end{array}\end{equation}
In particular, the bundle $\Lambda_\eta^1$ is defined by exactness of the
middle row and its sections are $1$-forms along the fibres of~$\eta$. Let us
define a vector bundle $A^1$ on $M$ by the exact sequence
\begin{equation}\label{defofA}
0\to\eta^*\Lambda_Z^{1,0}\to\Lambda_M^1\to A^1\to 0.\end{equation}
Then (\ref{diagram}) gives rise to two short exact sequences involving~$A^1$,
namely
\begin{equation}\label{firstAsequence}
0\to\eta^*\Lambda_Z^{0,1}\to A^1\to\Lambda_\eta^1\to 0\end{equation}
and
\begin{equation}\label{secondAsequence}
0\to B^1\to A^1\to\Lambda_M^{0,1}\to 0.\end{equation}
Observe that the short exact sequence (\ref{defofA}) together with the de~Rham
complex $\Lambda_M^\bullet$ on $M$ induces a differential complex $A^\bullet$
on~$M$. Each of (\ref{firstAsequence}) and (\ref{secondAsequence}) now gives
rise to a spectral sequence for computing the cohomology of $A^\bullet$ and
Theorem~\ref{mainthm} will be a consequence of combining these spectral
sequences with appropriate vanishing results, namely Jurchescu's
Theorem~\ref{vanishingthm} and a theorem of Buchdahl~\cite{b}.

The short exact sequence (\ref{secondAsequence}) induces a filtering of the
complex $A^\bullet$ on~$M$ and, in particular, we recover the resolution
(\ref{partialdolbeault}) and, indeed, coupled resolutions
\begin{equation}\label{coupledpartialdolbeault}
0\to\oe(B^p)\to{\mathcal E}_M^{0,0}(B^p)
\stackrel{\bar\partial}{\longrightarrow}
{\mathcal E}_M^{0,1}(B^p)\stackrel{\bar\partial}{\longrightarrow}
{\mathcal E}_M^{0,2}(B^p)\stackrel{\bar\partial}{\longrightarrow}\cdots
\stackrel{\bar\partial}{\longrightarrow}{\mathcal E}_M^{0,n}(B^p)\to 0,
\end{equation}
for each $p\geq 0$. In particular, we see that $B^p$ is partially holomorphic,
as presaged just prior to Theorem~\ref{mainthm}. The spectral sequence
associated to this filtration of $A^\bullet$ is
\begin{equation}\label{secondSS}
E_1^{p,q}=H^q(M,\oe(B^p))\Longrightarrow H^{p+q}(\Gamma(M,A^\bullet)).
\end{equation}
As with Stein manifolds, Theorem~\ref{vanishingthm} applies equally to the
cohomology of partially holomorphic vector bundles on a Cartan manifold. Hence,
the spectral sequence (\ref{secondSS}) collapses to an isomorphism
\begin{equation}\label{secondisomorphism}
H^p(\Gamma(M,A^\bullet))\cong H^p(\Gamma(M,\oe(B^\bullet))).\end{equation}

The short exact sequence (\ref{firstAsequence}) induces a different filtering
of the complex $A^\bullet$ on~$M$. In particular, the induced complexes
${\mathcal E}_\eta^\bullet(\eta^*\Lambda^{0,p})$ resolve the sheaves
$\eta^{-1}{\mathcal E}^{0,p}$ of germs of smooth sections of
$\eta^*\Lambda^{0,p}$ locally pulled back from~$Z$. The differential operator
$d_\eta:\eta^*\Lambda^{0,p}\to\Lambda_\eta^1(\eta^*\Lambda^{0,p})$ is sometimes
called a Bott partial connection. The spectral sequence associated to this
alternative filtration of $A^\bullet$ is
\begin{equation}\label{firstSS}
E_1^{p,q}=H^q(M,\eta^{-1}{\mathcal E}^{0,p})\Longrightarrow
H^{p+q}(\Gamma(M,A^\bullet)).
\end{equation}
In~\cite{b}, Buchdahl shows that if $\eta$ has connected fibres whose
de~Rham cohomology $H^q(\eta^{-1}(z),{\mathbb C})$ vanishes
for $q\geq 1$, then $H^q(M,\eta^{-1}{\mathcal E}^{0,p})=0$ for $q\geq 1$.
Hence, the spectral sequence (\ref{firstSS}) collapses to an isomorphism
\begin{equation}\label{firstisomorphism}
H^p(\Gamma(M,\eta^{-1}{\mathcal E}^{0,\bullet}))
\cong H^p(\Gamma(M,A^\bullet)).\end{equation}
Since $\eta$ has connected fibres, sections of $\eta^*\Lambda^{0,p}$ that are
locally pulled back from $Z$ are, in fact, pulled back from $Z$. Hence, the
left hand side of (\ref{firstisomorphism}) is the Dolbeault cohomology of~$Z$.
Combining (\ref{firstisomorphism}) and (\ref{secondisomorphism}) gives
$$H^p(Z,{\mathcal O})\cong H^p(\Gamma(M,\eta^{-1}{\mathcal E}^{0,\bullet}))
\cong H^p(\Gamma(M,A^\bullet))\cong
H^p(\Gamma(M,\oe(B^\bullet))),\ \forall p,$$
as required. \end{proof}

It is worth recording separately various special cases. First of all, there is
the spectral sequence (\ref{secondSS}) combined with the
isomorphism~(\ref{firstisomorphism}) when $M$ is a general mixed manifold:--
\begin{thm}\label{specseq} Suppose that the fibres of $\eta$ are contractible.
Then, there is a spectral sequence
$$E_1^{p,q}=H^q(M,\oe(B^p))\Longrightarrow H^{p+q}(Z,{\mathcal O}).$$
\end{thm}
\noindent There are also three interesting special cases of the theorem
itself:--
\begin{itemize}
\item we take $M$ to be $Z$ but regarded as a smooth manifold: in this case
$\oe(B^\bullet)$ reduces to the usual Dolbeault resolution;
\item we take $M$ to be a Stein manifold: in this case we obtain
the holomorphic language of~\cite{egw,egw1,g1};
\item we take $M\subset\Xi\times Z$ as in \S\ref{simple}, supposing that $M$ is a
Cartan manifold: we obtain the smoothly parameterised \v{C}ech cohomology and
Theorem~\ref{SmoothCechTheorem}.
\end{itemize}

As in~\S\ref{simple}, it is possible to be more explicit concerning the
isomorphism (\ref{mainresult}) in the presence of a section:--
\begin{thm}\label{sectionthm} Suppose $\gamma:Z\to M$ is a smooth section of
the submersion $\eta:M\to Z$. Then there is a mapping of complexes
\begin{equation}\label{chain}
\Gamma(M,\oe(B^\bullet))\longrightarrow\Gamma(Z,{\mathcal E}^{0,\bullet})
\end{equation}
that, if the hypotheses of Theorem~\ref{mainthm} are satisfied, induces the
isomorphism {\rm (\ref{mainresult})} on cohomology.
\end{thm}
\begin{proof} Our construction follows (\ref{pullback}) except that more care
has to be exercised, taking into account that $M$ is no longer locally a
product. {From} (\ref{secondAsequence}) and (\ref{defofA}) we see that
$$\omega\in\Gamma(M,B^p)\hookrightarrow\Gamma(M,A^p)$$
may be lifted to a $p$-form $\widetilde\omega$ on $M$ and then its pullback
$\gamma^*\widetilde\omega$ is a $p$-form on~$Z$. {From} (\ref{defofA}) it
follows that the $(0,p)$-component $(\gamma^*\widetilde\omega)^{0,p}$ is
independent of choice of lift. By construction, (\ref{chain}) is a composition
$$\Gamma(M,\oe(B^\bullet))\longrightarrow\Gamma(M,A^\bullet)
\longrightarrow\Gamma(Z,{\mathcal E}^{0,\bullet})$$
each of is evidently a mapping of complexes. In fact, the first of these has
nothing to do with $\gamma$ and, in the proof of Theorem~\ref{mainthm}, induces
the isomorphism~(\ref{secondisomorphism}). On the other hand,
(\ref{firstAsequence}) induces a chain mapping
$\eta^*:\Gamma(Z,{\mathcal E}^{0,\bullet})\to\Gamma(M,A^\bullet)$, which is the
source of the isomorphism with Dolbeault cohomology when the appropriate
topological conditions on the fibres of $\eta$ are satisfied. Since
$\eta\circ\gamma={\mathrm{Id}}$, we see that
$\gamma^*\circ\eta^*={\mathrm{Id}}$ and, in this case, taking the
$(0,p)$-component is evidently vacuous. \end{proof}

\section{Summary}\label{summary}
For those readers concerned only with the case of smoothly parameterised
\v{C}ech cohomology, it is perhaps worthwhile to summarise the events of
\S\S\ref{mixed}--\ref{general} as they apply to this case. In effect, the proof
boils down to the following. Let $\eta^{-1}{\mathcal O}$ denote the inverse
image sheaf of~${\mathcal O}$. It is the sheaf of smooth functions on $M$ locally
of the form $f\circ\eta$ for $f$ a holomorphic function on~$Z$. Equivalently,
they are locally constant along the fibres of $f$ or, in other words,
annihilated by~$d_\eta$, the exterior derivative along the fibres of~$\eta$. We
obtain a resolution
$$0\to\eta^{-1}{\mathcal O}\to
\oe(B^0) \stackrel{d_\eta}{\longrightarrow}\oe(B^1)
\stackrel{d_\eta}{\longrightarrow}\oe(B^2)\to\cdots$$
of $\eta^{-1}{\mathcal O}$ by sheaves that, by dint of Jurchescu's
Theorem~\ref{vanishingthm}, are acyclic when $M$ is Cartan. We conclude that
$$\frac{\ker:\Gamma(M,\oe(B^p))
             \stackrel{d_\eta}{\longrightarrow}\Gamma(M,\oe(B^{p+1}))}
       {\mbox{im : }\Gamma(M,\oe(B^{p-1}))
             \stackrel{d_\eta}{\longrightarrow}\Gamma(M,\oe(B^{p}))}
\cong H^p(M,\eta^{-1}{\mathcal O}).$$
On the other hand, there is always a tautological homomorphism
$$H^p(Z,{\mathcal O})\to H^p(M,\eta^{-1}{\mathcal O}),$$
which, by Buchdahl's theorem~\cite{b}, is an isomorphism when $\eta$ has
contractible fibres.

\section{Non-convex tube domains}\label{tubes}
In this section we discuss the application of Theorem~\ref{mainthm} to describe
the cohomology of tubes over non-convex cones, cf.~\cite{g}. Let
$V\subset{\mathbb R}^n$ be an open cone (so $x\in V$ $\implies$
$\lambda x\in V$ for $\lambda>0$). The $V$ we have in mind need not be convex.
In contrast to convex cones, non-convex cones can display very pathological
behaviour and we need to restrict the class of cones under consideration. We
fix a natural number $s$ and consider cones which are unions of $s$-planes.
More specifically, we shall employ smoothly parameterised Stein coverings
obtained as tubes over some special convex subcones of $V$ as follows. Let us
call a cone $W$ an $s$-wedge if it is the direct product of an $s$-subspace
$d(W)\in{\mathrm{Gr}}_s({\mathbb R}^n)$ and an $(n-s)$-dimensional sharp convex
cone. More precisely, $C$ should be a convex cone that does not contain a line
in an $(n-s)$-dimensional linear subspace complementary to $d(W)$ and $W$
should consist of those vectors in ${\mathbb{R}}^n$ whose projection onto this
complementary subspace lie in~$C$, in which case we shall write
$W=d(W)\times C$. Let us say that a cone $V$ is $s$-convex if it is the union
of $s$-wedges $W_\xi$ for $\xi\in\Xi$ some connected smooth manifold and
\begin{itemize}
\item the mapping $d:\Xi\to{\mathrm{Gr}}_s({\mathbb R}^n)$ is a finite-to-one
smooth covering with connected open range;
\item for each $x\in V$, the set $\{\xi\in\Xi\,\mbox{ s.t. }W_\xi\ni x\}$ is
non-empty contractible.
\end{itemize}

It is natural to expect for $s$-convex tubes that the interesting cohomology
occurs in degree~$s$. The notion of $0$-convex coincides with convex and,
certainly, it is natural to consider holomorphic functions since they provide
an infinite-dimensional function space whilst higher cohomology vanishes. Under
suitable conditions, it is possible to develop a very advanced theory of the
cohomology $H^s(T,{\mathcal{O}})$ of the corresponding tube domain
$Z=V+i{\mathbb{R}}^n$ analogous to Bochner's theory in convex tubes~\cite{g3}.

\medbreak\noindent{\bf Definition}\quad {\em If $V\subset{\mathbb{R}}^n$ is an
$s$-convex cone, we shall say that\:\ $V+i{\mathbb R}^n\subset{\mathbb C}^n$ is
an \underbar{$s$-convex tube}.}

\begin{thm} Suppose $Z=V+i{\mathbb R}^n$ is an $s$-convex tube. Define
$M\subset\Xi\times Z$ by
$$M=\{(\xi,z)\in\Xi\times Z\,\mbox{\rm\ s.t. }
z\in W_\xi\times i{\mathbb R}^n\}.$$
Then $M$ is a Cartan manifold and the natural projection $\eta:M\to Z$ has
contractible fibres.
\end{thm}
\begin{proof} The fibre of the natural projection $\tau:M\to\Xi$ over $\xi$ is
$W_\xi+i{\mathbb R}^n$. This is convex, so Proposition~\ref{supply}
implies $M$ is Cartan. The fibre of $\eta$ over $z=x+iy\in Z$ is
$\{\xi\in\Xi\,\mbox{ s.t. }W_\xi\ni x\}$. This is contractible by definition of
$s$-convex. \end{proof}
Immediately from Theorem~\ref{SmoothCechTheorem}, we conclude:--
\begin{cor} The analytic cohomology of an $s$-convex tube may be computed by
smooth \v{C}ech cohomology.
\end{cor}
Examples of $s$-convex tubes are provided by pseudo-Hermitian symmetric
spaces of tube type~\cite{fg}. Here we describe three such spaces:--

\subsection{First example}\label{firsteg}
Generalising \S\ref{san}, let $n=p+q$ with $p\geq 2$, introduce on
${\mathbb{R}}^n$ the non-degenerate non-degenerate symmetric form
$\langle\phantom{x},\phantom{x}\rangle$ of type $(p,q)$
$$\langle x,y\rangle=x_1y_1+\cdots+x_py_p
-x_{p+1}y_{p+1}-\cdots-x_ny_n,$$
and set
$$V=\left\{x\in{\mathbb{R}}^n\mbox{ s.t. }\langle x,x\rangle>0\right\}.$$
This is a connected non-convex cone and we maintain that is $(p-1)$-convex, as
follows.
Choose an orientation on one, and hence all, $p$-dimensional linear subspaces
on which $\langle\phantom{x},\phantom{x}\rangle$ is positive definite. Let
$\xi$ denote an oriented $(p-1)$-dimensional linear subspace of
${\mathbb{R}}^n$ on which $\langle\phantom{x},\phantom{x}\rangle$ is positive
definite. Then $\langle\phantom{x},\phantom{x}\rangle$ restricted to
$\xi^\perp$ is a non-degenerate symmetric form of type $(1,q)$. Therefore,
$$C_\xi=\{x\in{\mathbb{R}}^n\mbox{ s.t. }
\langle x,\xi\rangle=0,\;\langle x,x\rangle>0,\;\mbox{and
$\xi\oplus{\mathbb{R}}x$ has chosen orientation}\}$$
is a sharp convex cone, whence $W_\xi=\xi\times C_\xi$ is a $(p-1)$-wedge.
In this way, V is covered by $(p-1)$-wedges parameterised by
$$\Xi=\left\{\xi\in{\mathrm{Gr}}_{p-1}^+({\mathbb{R}}^n)\mbox{ s.t. }
\langle\phantom{x},\phantom{x}\rangle|_\xi\mbox{ is positive definite}
\right\}$$
where ${\mathrm{Gr}}^+$ denotes the oriented Grassmannian. The mapping
$d:\Xi\to{\mathrm{Gr}}_{p-1}({\mathbb{R}}^n)$ simply forgets the orientation.
Notice that the corresponding smoothly parameterised \v{C}ech covering of the
$(s-1)$-convex tube $V+i{\mathbb{R}}^n$ respects the action
of the natural symmetry group ${\mathrm{SO}}(p,q)$.

\subsection{Second example}
Let $n=m(m+1)/2$ and identify ${\mathbb R}^n$ with the $m\times m$ real
symmetric matrices. Let $m=p+q$ and set
$$V=\{X\in{\mathbb R}^n\mbox{ s.t. $X$ has signature }(p,q)\}.$$
Then $Z=V+i{\mathbb R}^n$ is a $pq$-convex tube. Matrices of the form
$$X=\left\lgroup\begin{array}{cc}0&B\\ B^t&0\end{array}\right\rgroup+
\left\lgroup\begin{array}{cc}P&0\\ 0&Q\end{array}\right\rgroup,$$
where $B$ is an arbitrary $p\times q$ matrix, $P$ is $p\times p$ symmetric
positive definite, and $Q$ is $q\times q$ symmetric negative definite, form a
$pq$-wedge and a covering by $pq$-wedges may be obtained by moving this one
under the action $X\mapsto AXA^t$ of ${\mathrm{SL}}(m,{\mathbb{R}})$.

\subsection{Third example}
Let $n=m^2$, identify ${\mathbb R}^n$ with the $m\times m$ real
matrices, and set
$$V=\{X\in{\mathbb R}^n\mbox{ s.t. }\det X >0\}.$$
Then $Z=V+i{\mathbb R}^n$ is an $m(m-1)/2$-convex tube. There is a  basic
wedge in~$V$:--
$$W=\{X=S+P\mbox{ s.t. $S$ is skew and $P$ is symmetric positive definite}\}$$
and the general wedge in the smoothly parameterised covering is obtained from
this one under the action $X\mapsto AXA^{-1}$ of
${\mathrm{SL}}(m,{\mathbb{R}})$.

\subsection{Example \S\ref{san} revisited}
Whilst the tube over the non-convex cone
$$V=\{(x_1,x_2,x_3)\in{\mathbb{R}}^3\mbox{ s.t. }x_1{}^2+x_2{}^2>x_3{}^2\}$$
admits a smoothly parameterised Stein covering, it also admits a foliation by
closed Stein submanifolds as follows. Firstly, define $\pi:V\to S^1$ by
$$\begin{array}{rcl}(x_1,x_2,x_3)&\stackrel{\pi}{\longrightarrow}&\displaystyle
\left(\frac{x_1x_3+x_2\sqrt{x_1{}^2+x_2{}^2-x_3{}^2}}{x_1{}^2+x_2{}^2},
\frac{x_2x_3-x_1\sqrt{x_1{}^2+x_2{}^2-x_3{}^2}}{x_1{}^2+x_2{}^2}\right)\\[12pt]
&=&(p,q)\quad\mbox{characterised by }\left\{ \begin{array}lp^2+q^2=1\\
px_1+qx_2=x_3\\
px_2-qx_1>0.
\end{array}\right.\end{array}$$
For each $(p,q)\in S^1$, the inverse image $\pi^{-1}(p,q)$ is a half plane
tangent to the null cone $\{x_1{}^2+x_2{}^2=x_3{}^2\}$. Now define
$\tau:V+i{\mathbb{R}}^3=Z\to\Xi=S^1\times{\mathbb{R}}$ by
$$\tau(x_1+iy_1,x_2+iy_2,x_3+iy_3)=(p,q,py_1+qy_2-y_3),\quad
\mbox{where }(p,q)=\pi(x_1,x_2,x_3).$$
Then, $\tau^{-1}(p,q,0)$ is a half-plane in the complex linear subspace
$\{pz_1+qz_2=z_3\}$ and, more generally, $\tau^{-1}(p,q,r)$ is an imaginary
translate thereof. In particular, the fibres of $\tau$ are all Stein. The
submersion $\tau:Z\to\Xi$ evidently endows $Z$ with the structure of a mixed
manifold of type $(2,2)$. Let us write $M$ for $Z$ endowed with this structure.
Arguing as in the proof of Proposition~\ref{supply}, we see that $M$ is a
Cartan manifold. If we denote by $\eta:M\to Z$ the identity mapping, then we
may use Theorem~\ref{mainthm} to realise $H^p(Z,{\mathcal{O}})$ as follows.
Firstly, let us write $\Lambda_\tau^{0,1}$ instead of $\Lambda_M^{0,1}$ for the
$(0,1)$-forms along the fibres of~$\tau$. Then, because the bundle
$\Lambda_\eta^1$ vanishes, the two short exact sequences (\ref{firstAsequence})
and (\ref{secondAsequence}) reduce to the single short exact sequence
$$0\to B^1\to\Lambda_Z^{0,1}\to\Lambda_\tau^{0,1}\to 0.$$
In our case, $B^1$ is a line bundle. The $\bar\partial$-operator on
$Z$ induces
$\bar\partial_\tau:B^1\to\Lambda_\tau^{0,1}\otimes B^1$
defining the partially holomorphic structure on~$B^1$. Also $\bar\partial$
induces $d_\eta:\oe\to\oe(B^1)$ and Theorem~\ref{mainthm} says that
$\Gamma(Z,{\mathcal{O}})=\ker d_\eta:\Gamma(Z,\oe)\to\Gamma(Z,\oe(B^1))$ and
$$H^1(Z,{\mathcal{O}})=
{\mathrm{coker}}\,d_\eta:\Gamma(Z,\oe)\to\Gamma(Z,\oe(B^1)).$$
This last statement implies that we may always find a Dolbeault representative
$\omega$ on $Z$ that restricts to zero on every leaf of our Stein foliation.
Another immediate consequence of Theorem~\ref{mainthm} is that all higher
cohomology vanishes.

\section{Twistor spaces}\label{twistors}
In order to have a true example of Theorem~\ref{mainthm}, other than smoothly
parameterised \v{C}ech cohomology, in this section we sketch how the cohomology
of the twistor space of a self-dual Riemannian manifold may be described and
how this description may be used to invert the Penrose transform. The twistor
space $Z$ of a self-dual Riemannian manifold $\Xi$ is constructed
in~\cite{ahs}. Its salient features are as follows. The sphere bundle in the
anti-self-dual $2$-forms is naturally a $3$-dimensional complex manifold $Z$
and the fibres of the mapping $\pi:Z\to\Xi$ are~${\mathbb{CP}}_1$'s. There is
an anti-holomorphic involution $\sigma:Z\to Z$, which is the antipodal map on
each fibre of~$\pi$. There is a holomorphic line-bundle ${\mathcal{O}}(1)$ on
$Z$, which is the hyperplane section bundle on each fibre of~$\pi$. The Penrose
transform~\cite{epw,h,w} identifies the cohomology $H^1(Z,{\mathcal{O}}(k))$
with solutions of certain conformally invariant systems of partial differential
equations on~$\Xi$, the so-called `massless field equations'. In \cite{w}
Woodhouse establishes this result by utilising canonical representatives of the
Dolbeault cohomology, characterised as being harmonic on each fibre of~$\pi$.
In fact, in the conformally flat case, these canonical representatives coincide
with those already constructed by Gindikin and Henkin via complexification and
the classical $\kappa$-operator from integral geometry
in~\cite{gh1,gh2,gh3,gh4}.

Define
$M\equiv\{(z,w)\in Z\times Z\mbox{ s.t. }\pi(z)=\pi(w)\mbox{ but }z\not= w\}$.
Then, we have the commutative diagram:--
\begin{center}\begin{picture}(80,50)
\put(0,5){\makebox(0,0){$Z$}} \put(40,45){\makebox(0,0){$M$}}
\put(80,5){\makebox(0,0){$\Xi$}} \put(30,35){\vector(-1,-1){20}}
\put(50,35){\vector(1,-1){20}} \put(15,30){\makebox(0,0){$\eta$}}
\put(65,30){\makebox(0,0){$\tau$}}
\put(10,5){\vector(1,0){60}}
\put(40,10){\makebox(0,0){$\pi$}}
\end{picture}\end{center}
where $\eta(z,w)=z$. The fibres of $\tau$ are intrinsically
${\mathbb{CP}}_1\times{\mathbb{CP}}_1\setminus\Delta$ where $\Delta$ is the
diagonal. In particular, they are holomorphic and, indeed, Stein. In fact, $M$
is a Cartan manifold. The mapping $\eta:M\to Z$ is a submersion with
contractible fibres (intrinsically each fibre
is~${\mathbb{CP}}_1\setminus\{\mbox{point}\}$). Therefore, we are in a position
to use Theorem~\ref{mainthm}. {From} a massless field on $\Xi$, it is
straightforward to identify the bundles $B^p$ on $M$ and, at least for
right-handed fields, write down a specific representative in the cohomology
$H^1(\Gamma(M,\oe(B^\bullet)))$ inverting the Penrose transform. Pulling back
this representative under $\gamma:Z\to M$ given by $\gamma(z)=(z,\sigma(z))$
and taking the $(0,1)$-component in accordance with Theorem~\ref{sectionthm}
gives the Woodhouse representative~\cite{w}. {From} the point of view of this
article, Theorem~\ref{mainthm} inverts the Penrose transform in the Riemannian
setting (directly going between cohomology on $Z$ and fields on $\Xi$), whereas
the holomorphically parameterised \v{C}ech cohomology of \cite{egw,egw1,g1}
(which Theorem~\ref{mainthm} generalises) inverts the Penrose
transform via the complexification of~$\Xi$, utilising that the massless field
equations on $\Xi$ are elliptic and whose solutions are thereby real-analytic.
This is the inversion via the $\kappa$-operator detailed
in~\cite{gh1,gh2,gh3,gh4}.

\end{document}